\documentclass{article}

\usepackage{amssymb,amscd,amsthm,verbatim,fancyhdr,mathrsfs}
\usepackage[leqno]{amsmath}
\usepackage{mathrsfs}

\begin{document}
\numberwithin{equation}{section}
\newtheorem{thm}{Theorem}[section]
\newtheorem{lemma}[thm]{Lemma}
\newtheorem{clm}[thm]{Claim}
\newtheorem{remark}[thm]{Remark}
\newtheorem{definition}[thm]{Definition}
\newtheorem{cor}[thm]{Corollary}
\newtheorem{prop}[thm]{Proposition}
\newcommand{\hdt}{{\dot{\mathrm{H}}^{1/2}}}
\newcommand{\hdtr}{{\dot{\mathrm{H}}^{1/2}(\mathbb{R}^3)}}
\newcommand{\R}{\mathbb{R}}
\newcommand{\rd}{{\mathbb{R}^d}}
\newcommand{\ei}{\mathrm{e}^{it\Delta}}
\newcommand{\ltrt}{L^3(\mathbb{R}^3)}
\newcommand{\lt}{{L^3}}
\newcommand{\lprd}{{L^p(\mathbb{R}^d)}}
\newcommand{\lp}{{L^p}}
\newcommand{\rt}{\mathbb{R}^3}
\newcommand{\X}{\mathfrak{X}}
\newcommand{\F}{\mathfrak{F}}
\newcommand{\hdhalf}{{\dot H^\frac{1}{2}}}
\newcommand{\hdthalf}{{\dot H^\frac{3}{2}}}
\newcommand{\hdo}{\dot H^1}
\newcommand{\rthmiz}{\R^3\times(-\infty,0)}
\newcommand{\q}[2]{{#1}_{#2}}
\renewcommand{\t}{\theta}
\newcommand{\lxt}[2]{L_{x,\,t}^{#1}}
\newcommand{\rr}{\sqrt{x_1^2+x_2^2}}
\newcommand{\ve}{\varepsilon}
\newcommand{\hdhrt}{\dot H^\frac{1}{2}(\mathbb{R}^3)}
\renewcommand{\P}{\mathbb{P} }
\newcommand{\RR}{\mathcal{R} }
\newcommand{\TT}{\overline{T} }
\newcommand{\e}{\epsilon }
\newcommand{\D}{\Delta }
\renewcommand{\d}{\delta }
\renewcommand{\l}{\lambda }
\newcommand{\To}{\TT_1 }
\newcommand{\ukt}{u^{(KT)} }
\newcommand{\ttil}{\tilde{T_1} }
\newcommand{\etl}{e^{t\Delta} }
\newcommand{\et}{\mathscr{E}_T }
\newcommand{\ft}{\mathscr{F}_T }
\newcommand{\eti}{\mathscr{E}^{\infty}_T }
\newcommand{\fti}{\mathscr{F}^{\infty}_T }
\newcommand{\xjn}{x_{j,n}}
\newcommand{\ljn}{\l_{j,n} }
\newcommand{\lkn}{\l_{k,n} }
\newcommand{\voj}{V_{0,j} }
\newcommand{\uon}{u_{0,n} }
\newcommand{\N}{\mathbb{N} }
\newcommand{\Z}{\mathbb{Z} }
\newcommand{\E}{\mathscr{E} }
\newcommand{\tu}{\tilde{u} }
\newcommand{\etj}{{\E_{T^*_j}}}
\newcommand{\B}{\mathcal{B}}
\newcommand{\C}{\mathbb{C}}
\newcommand{\tujn}{\tilde U_{j,n}}
\newcommand{\LL}{\Lambda}
\newcommand{\binf}{{\dot B^{-{d/p}}_{\infty,\infty}}}
\newcommand{\brq}{{\dot B^{s_{p,r}}_{r,q}}}
\newcommand{\bfin}{{\dot B^{s_{p,r}}_{r,q}}}
\newcommand{\tn}{{\tilde \|}}
\newcommand{\xto}{\xrightarrow[n\to\infty]{}}
\newcommand{\besa}{{\dot B^{s_{p,a}}_{a,q}}}
\newcommand{\besb}{{\dot B^{s_{p,b}}_{b,r}}}

\title{Profile decompositions for critical\\ Lebesgue and Besov space embeddings}
\author{{\sc Gabriel S. Koch}\\ \small Oxford University\\ \small {\em koch@maths.ox.ac.uk}}
\maketitle
\begin{abstract}
Profile decompositions for ``critical" Sobolev-type embeddings are established, allowing one to regain some compactness despite the non-compact nature of the embeddings.  Such decompositions have wide applications to the regularity theory of nonlinear partial differential equations, and have typically been established for spaces with Hilbert structure.  Following the method of S. Jaffard, we treat settings of spaces with only Banach structure by use of wavelet bases.  This has particular applications to the regularity theory of the Navier-Stokes equations, where many natural settings are non-Hilbertian.
\end{abstract}
\ \\

In this article, we characterize the causes for the lack of compactness in certain embeddings of Banach spaces.  In particular, one can re-write a bounded sequence in such a way that it is obvious why there fails to be a convergent subsequence (even in the weaker space).  Once the defects of compactness are thus identified, it is often possible to isolate them and regain some aspects of compactness.  \\

This useful tool for re-writing bounded sequences is known as a ``decomposition into profiles" (typically after passing to a subsequence), the ``profiles" being the typical obstacle to compactness -- specifically, norm-invariant transformations of fixed non-zero elements of the space.  Since these elements are fixed (for the whole sequence) one can think of them as ``limits" which can replace the need for an actual convergent subsequence.  This is particularly useful in the study of nonlinear partial differential equations, where the natural setting is often that of such Banach spaces, and the convergence of certain sequences corresponds to establishing the existence of certain solutions to the equation being considered.  In fact, it is primarily for this purpose that such decompositions have historically been established.\\

As we have just intimated, such a program is not new.  Patrick G\' erard \cite{gerard} was probably the first to do this and treated the embedding $\dot H^s(\rd) \hookrightarrow L^q(\rd)$ for $s>0$ and an appropriate $q>1$.  It is typical that such decompositions are established for spaces of functions defined on $\rd$, since the invariance of $\rd$ under translations and dilations leads to the typical defect of compactness of the embedding (see below). Later, St\' ephane Jaffard \cite{jaffard} extended the work of Gerard (although with slightly less specific results) to $\dot H^{s,p}(\rd) \hookrightarrow L^q(\rd)$ for $s>0$ and appropriate $p,q >1$.  His proof was somewhat technically different from \cite{gerard}, and he was able to move to the non-Hilbertian setting (i.e., spaces not based on $L^2$) by use of wavelet bases.  There have of course been many other works on profile decompositions, see for example \cite{BG1,bmm,bc,km1,keraani,keraani2,mv}, to name only a few. \\

In this article we follow the method of Jaffard to prove similar decompositions for the embedding $L^p \hookrightarrow \bfin$ for $p\geq 2$, where $\bfin$ is any one of a family of Besov spaces for $r,q>p$ with a particular negative index of regularity $s_{p,r}$.  Of course, in general these spaces are non-Hilbertian and are therefore well suited to the method described in \cite{jaffard}.  As an indication of the generality of the method (and with specific applications in mind which we'll describe momentarily), we show how similar results may be obtained as well for embeddings of the type $\besa \hookrightarrow \besb$. We also clarify, correct and improve some of the arguments in \cite{jaffard}, for example showing that a compactly supported wavelet basis is not needed and more carefully detailing the relationships between the ``orthogonal cubes" in the construction.\\

The method of establishing profile decompositions using wavelet bases is actually quite general and can be applied to many embeddings between spaces of Besov or Triebel-Lizorkin type, which we plan to present in a unified and hopefully simpler way in the upcoming work \cite{ck}.  Our motivation for studying the particular embeddings in this work now are their applications to the regularity theory of the Navier-Stokes equations, which we'll now proceed to describe.\\

In \cite{kk}, regularity results for Navier-Stokes were established based on the profile decomposition of Isabelle Gallagher \cite{gallagher} (based in turn on \cite{gerard}) in the setting of $\hdhalf(\rt)$.  Perhaps a more natural setting for these equations is the space $L^d(\rd)$, however in that work (and others, see e.g. \cite{sverakrusin}) the more special setting of $\hdhalf$ was investigated.  This is partly because much theory (such as \cite{gerard}) has been developed in the $L^2$ setting, and partly because the Hilbert structure of $L^2$-based spaces makes the analysis somewhat simpler.    In the non-Hilbertian setting, somewhat different methods are in general required. The results and methods presented below will make it possible to move to the more natural (non-Hilbertian) settings in the Navier-Stokes framework (see, e.g., the upcoming work \cite{gkp} where we shall give extensions of \cite{kk} and \cite{sverakrusin} in the more natural setting), and hopefully in the study of other PDEs as well for which $H^s$ is not the only interesting or natural setting.
\\

Navier-Stokes is of course not the only application of interest.  As we mentioned above, using profile decompositions has become quite a common approach in the study of many nonlinear partial differential equations in ``critical" function-space settings, see e.g. \cite{ckm,km1,km2,km3,kk,ks,mxz3} to name only a few, where the existence of a profile decomposition implies various regularity results for solutions of the PDE.  Until the recent work \cite{kk}, the methods had been applied exclusively to dispersive and hyperbolic PDE where the Hilbertian $H^s$ settings were the natural ones.  It was therefore natural to perform the analysis of the parabolic Navier-Stokes equations in \cite{kk} in a similar setting, namely $\hdhalf(\rt)$, and moreover until now the more general $L^p$ theory was not available.  The application of the results in this paper to Navier-Stokes in the natural setting of $L^d(\rd)$ will be addressed soon in \cite{gkp}, in the same spirit as \cite{kk}.  This will give a different proof of the well-known and difficult result in \cite{ess} (see also \cite{dongdu}), that the $L^3(\rt)$ norm of a solution to Navier-Stokes must become infinite near a singularity (the result of \cite{kk} treated the special case of blow-up of the stronger $\hdhalf$ norm).
\\

It is interesting to note as well that profile decompositions (which we will describe below) can also be used (see \cite{gerard}) to recover the existence of minimizers in the classical Sobolev embeddings, a fact established by P. L. Lions \cite{lions3} in the '80s using his method of ``concentration-compactness".
\\

Let us finally outline now what a profile decomposition is.  The principle defect in the lack of compactness of many embeddings  $X \hookrightarrow Y$ of homogeneous spaces of functions defined on $\rd$ comes from the following example:  For some non-zero element $f\in X$ and $\alpha >0$, for any $n\in \N$ and $x\in \rd$ define
$$f_n(x):=\frac{1}{(\l_n)^\alpha}f \left(\frac{x-x_n}{\l_n}\right)\ ,$$
where $\l_n >0$ and $x_n \in \rd$.  The spaces and the scaling constant $\alpha$ will be such that $\|f_n\|_X \equiv \|f\|_X$ and $\|f_n\|_Y \equiv \|f\|_Y$ (when the scalings of the spaces coincide in this way we call the embedding ``critical" -- in the context of PDE this usually means moreover that the scaling is a natural one for the equation).  Moreover, $X$ and $Y$ are such that if either $\lim_{n\to\infty}\l_n \in \{0,\infty\}$ or $\l_n \equiv 1$ and $\lim_{n\to\infty}|x_n| = \infty$, then $f_n \rightharpoonup 0$ in $X$ and $Y$.  Then one has $\sup_{n\in\N}\|f_n\|_X = \|f\|_X < \infty$, yet clearly one cannot have any subsequence which is strongly convergent (necessarily to zero) in $Y$ since $\|f_n\|_Y \equiv \|f\|_Y >0$.  A profile decomposition for bounded sequences in $X$ is a way to show that the lack of compactness of such embeddings comes in general from linear combinations of the above type of ``profile" example.  Moreover, the profiles do not interact with each other in ways which might ``cancel" the defect of compactness caused by each one.
\\
\pagebreak

The main example we study here is the following: Suppose $2\leq p < q,r \leq \infty$ and set $s_{p,r}:= d(\frac{1}{r} - \frac{1}{p}) < 0$.  Then $\lprd \hookrightarrow \brq$ and we have the following profile decomposition for bounded sequences in $\lprd$:
\begin{thm}\label{thm:a}
Let $\{u_n\}_{n=1}^\infty$ be a bounded sequence in $\lprd$.  There exists a sequence of profiles $\{\phi_l\}_{l=1}^\infty \subset \lprd$ and for each $l\in \N$ a sequence $\{(j_n^l,k_n^l)\}_{n=1}^\infty \subset \Z \times \Z^d$, both depending on $\{u_n\}$, such that, after possibly passing to a subsequence in $n$,
\begin{equation}\label{profiles}
u_n(x) = \sum_{l=1}^L (2^{j^l_n})^{{d/p}} \phi_l(2^{j^l_n}x - k^l_n) + r_n^L(x)
\end{equation}
for any $L \in \N$ where the following properties hold:
\begin{enumerate}
\item For $l\neq l'$, the sequences $\{(j_n^l,k_n^l)\}$ and $\{(j_n^{l'},k_n^{l'})\}$ are orthogonal in the following sense:
\begin{equation}\label{orth1}
\left|\log \left(2^{(j_n^l - j_n^{l'})}\right)\right| + \left|2^{(j_n^l - j_n^{l'})}k_n^{l'} - k_n^l\right|  \xrightarrow[n\to\infty]{} +\infty
\end{equation}
\item The remainder $r_n^L$ satisfies the following smallness condition:
\begin{equation}\label{orth2}
\lim_{L\to\infty} \left(\varlimsup_{n\to \infty} \|r_n^L\|_\brq \right) = 0
\end{equation}
\item There is a norm $\| \cdot \tilde \|_\lprd$ which is equivalent to $\|\cdot \|_\lprd$ such that for each $n \in \N$,
\begin{equation}\label{orth3}
\sum_{l=1}^\infty \|(2^{j^l_n})^{{d/p}} \phi_l(2^{j^l_n}x - k^l_n) \tn_\lprd^p \leq \varliminf_{n'\to\infty} \|u_{n'}\tn_\lprd^p
\end{equation}
and for any $L \in \N$,
\begin{equation}\label{orth4}
 \|r_n^L\tn_\lprd \leq \|u_n\tn_\lprd + \circ(1) \quad \textrm{as} \quad n\to\infty\ .
\end{equation}
\end{enumerate}
\end{thm}

The orthogonality property (\ref{orth1}) implies exactly that the interaction between the profiles $\{\phi_l\}$ (after the norm-invariant transformations) becomes negligible as $n\to\infty$ (see e.g. \cite{gallagher,kk} as well as discussions below).
Properties (\ref{orth3}) and (\ref{orth4}) can be thought of as ``stability properties" of the profile decomposition in terms of the sequences from which they are derived.  We also comment here that the norm introduced in (\ref{orth3}) - (\ref{orth4}) may not have the same transformational invariance as $\| \cdot \|_\lprd$, but (\ref{orth3}) does imply that there exists some universal constant $C>0$ such that
\begin{equation}\label{altorth}
\sum_{l=1}^\infty \|\phi_l \|_\lprd^p \leq C \varliminf_{n\to\infty} \|u_{n}\|_\lprd^p
\end{equation}
due to the invariance of the usual $\lprd$ norm under the transformations on $\phi_l$ in (\ref{orth3}).  One could alternatively define a new $\tilde{\phi}_l:=(2^{j^l_1})^{{d/p}} \phi_l(2^{j^l_1}\cdot - k^l_1)$ (and adjust $\{(j_n^l,k_n^l)\}$ in the obvious way) to obtain (\ref{altorth}) with the norm $\| \cdot \tn_\lprd$ and $C=1$ (plus all the previous properties), which gives a more standard and useful property (see e.g. \cite{gkp}).  With minimal adjustments, a similar theorem about bounded sequences in Besov norms (which also has applications to Navier-Stokes, see \cite{gkp}) is also seen to be true, see Theorem \ref{thm:aa}.
\\

As mentioned above, the method of proof we will use follows that of \cite{jaffard}, using the ``wavelet" norms $\|\cdot \tn_\lprd$, $\| \cdot \tn_\brq$, etc., which are based on the expansion of functions in all spaces in {\em one} unconditional wavelet basis (and are equivalent to the standard norms).  Roughly, one extracts and categorizes the ``largest" wavelet components of $u_n$ in an exhaustive manner, which decreases the norms of $u_n$ in a ``largest" space $Z$ (consisting simply of functions with bounded wavelet coefficients) while the norms in the ``smallest" space $X$ (e.g., $\lprd$) remain bounded.  In fact, the remainder must tend to zero in $Z$ due to properties of ``nonlinear (wavelet) projections".  Then an interpolation inequality of the form $\|f\|_Y \lesssim \|f\|_X^\alpha \|f\|_Z^{1-\alpha}$ for some ``middle" space $Y$ (e.g. $Y=\brq$ and $X \hookrightarrow Y \hookrightarrow Z$) completes the main part of the proof.\\

{\bf Acknowledgements:}
\quad
The author would like to thank Professors Albert Cohen, Isabelle Gallagher and St\' ephane Jaffard for many helpful discussions.  \\

This work was supported by the EPSRC Science and Innovation award to the Oxford Centre for Nonlinear PDE (EP/E035027/1).
\section{Preliminaries}
\noindent
We'll start with the following facts regarding wavelet bases (see e.g. \cite{daubechies}, \cite{meyer2}; see also \cite{bf} and \cite{lr} for divergence-free wavelet bases):
\\\\
Fix any $m\in \N$.  Denote
$$\l = (i,j,k) \in \{1, \dots , 2^d -1\} \times \Z \times \Z^d =: \LL\ .$$
There exists a real-valued set of functions\footnote{In fact, one can take $\varphi \in \mathcal{C}^m_0$, but we do not presently require compact support of the basis functions.} $\{ \varphi^{(i)} \}_{1\leq i \leq 2^d -1} \subset \mathcal{C}^m (\R^d) \cup (\cap_{p>1}L^p(\R^d))$ such that $\varphi_\l = \varphi^{(i)}_{j,k}$ defined by $\varphi_\l(x) = 2^\frac{jd}{2} \varphi^{(i)}(2^j x -k)$ is an orthonormal basis of $L^2$:
$$f\in L^2 (\R^d) \iff f= \sum_{\l \in \LL} c_\l \varphi_\l, \qquad \textrm{where} \quad c_\l = c^{(i)}_{j,k} := \int_{\R^d} \varphi_\l \cdot f\ .$$
Moreover, we have the following equivalences of norms for other function spaces:
\begin{equation}\label{a}
\|f\|_{L^p(\R^d)} \simeq \|f\tilde\|_{L^p(\R^d)} := \left\| \left\{ \sum_{\l \in \LL} |c_{j,k}^{(i)}|^2 2^{dj} \chi_\d(x) \right\}^\frac{1}{2}\right\|_{L^p(\R^d)}
\end{equation}
for $1< p < \infty$ where $\chi_\d$ is the indicator function of the
cube\footnote{Note that these cubes correspond roughly to the ``supports" of the $\psi_\l$s.}
$\d = \d(\l)  =\d^{j,k}= {\{x \in \R^d \ | \ 2^j x -k \in [0,1)^d \}}$ (note that $|\d^{j,k}| = 2^{-dj}$), and
\begin{equation}\label{b}
\|f\|_{\dot B^s_{a,b}} \simeq \|f\tilde\|_{\dot B^s_{a,b}} := \left\| 2^{j(s + d(\frac{1}{2} - \frac{1}{a}))} \|c^{(i)}_{j,k}\|_{\ell^a_{i,k}} \right\|_{\ell^b_j}
\end{equation}
for $1 \leq a,b \leq \infty$ so long as $|s| < m$.  Here and in what follows we use the notation $A(f) \simeq B(f)$ whenever there exists $c>0$ independent of $f$ such that $c^{-1} A(f) \leq B(f) \leq c A(f)$ and $A(f) \lesssim B(f)$ if at least the first inequality holds.
\\\\
In particular, we have
$$\|f\tilde \|_{\dot B^{-{d/p}}_{\infty,\infty}} = \sup_{\l \in \LL}  2^{jd(\frac{1}{2} - \frac{1}{p})} |c^{(i)}_{j,k}| $$
when $m>d$, which motivates us to define $a^{(i)}_{j,k} := 2^{jd(\frac{1}{2} - \frac{1}{p})} c^{(i)}_{j,k}$ and $\psi_\l$ by $\psi_\l(x) = 2^{jd(\frac{1}{p} - \frac{1}{2})} \varphi_\l (x)$, so that for a given $f$ one may write $f=\sum_\l c_\l \varphi_\l = \sum_\l a_\l \psi_\l$.  With these definitions, we have the equivalence
\begin{equation}\label{c}
\|f\|_{\dot B^s_{a,b}} \simeq \|f\tilde\|_{\dot B^s_{a,b}} = \left\| 2^{j(s + d(\frac{1}{p} - \frac{1}{a}))} \|a^{(i)}_{j,k}\|_{\ell^a_{i,k}} \right\|_{\ell^b_j}
\end{equation}
so that in particular
\begin{equation}\label{d}
\|f\|_{\dot B^{-{d/p}}_{\infty,\infty}} \simeq \|f\tilde \|_{\dot B^{-{d/p}}_{\infty,\infty}}=  \sup_{\l \in \LL} |a^{(i)}_{j,k}|\ ,
\end{equation}
and
\begin{equation}\label{e}
\|f\|_{L^p(\R^d)} \simeq \|f\tilde\|_{L^p(\R^d)} = \left\| \left\{ \sum_{\l \in \LL} |a_{j,k}^{(i)}|^2 2^{\frac{2d}{p}j} \chi_\d(x) \right\}^\frac{1}{2}\right\|_{L^p(\R^d)}\ .
\end{equation}
Note that for any $(i,j,k) \in \LL$, by (\ref{e}) we have
\begin{equation}\label{f}
|a^{(i)}_{j,k}|^p =  \int_{\rd} \left\{ |a_{j,k}^{(i)}|^2 2^{\frac{2d}{p}j} \chi_\d(x) \right\}^\frac{p}{2}\ dx \leq \|f\tilde\|_\lprd^p\ .
\end{equation}
Hence $\lprd \hookrightarrow \dot B^{-{d/p}}_{\infty,\infty}$, and for $\{u_n\} \subset \lprd$ with $u_n = \sum_\l a_{\l,n}\psi_\l$ we see that
\begin{equation}\label{k}
\sup_n \|u_n\tilde \|_\lprd \leq C \quad \Longrightarrow \quad \sup_n \|u_n\tn_\binf = \sup_{\l,n}|a_{\l,n}| \leq C\ ,
\end{equation}
i.e. boundedness in $\lprd$ implies uniform boundedness of the wavelet coefficients by (\ref{d}).
\\\\
We'll need the following two facts regarding the wavelet coefficients of elements of $\lprd$:
\begin{prop}\label{prop1}
Let $\{a_\l\}$ correspond to some $f\in \lprd$.  Then there exists some $\l_0 \in \LL$ such that $\sup_{\l \in \LL} |a_\l| = |a_{\l_0}|$.
\end{prop}
\noindent
The proposition is a simple consequence of (\ref{e}), but for the convenience of the reader we give a proof in the Appendix, see Proposition \ref{prop:a1}.
\\\\
For the second fact, we will need the following definition of ``nonlinear projections":
\\\\
For any $u\in \binf$ and $N \in \N$, we'll say $P_N(u) \in \bar{P}_N(u)$ if
$P_N(u) = \sum_{m=1}^N a_{\l_m} \psi_{\l_m}$ for some ordering $\{\l_m\}_{m=1}^\infty$ of $\LL$ such that
\begin{equation}\label{g}
u = \sum_{m=1}^\infty a_{\l_m} \psi_{\l_m} \qquad \textrm{and} \qquad |a_{\l_m}| \geq |a_{\l_{m+1}}| \ \ \forall m\geq 1\ .
\end{equation}
Note that (\ref{g}) is always possible by Proposition \ref{prop1}.  We'll say that $Q_N(u) \in \bar{Q}_N(u)$ if $Q_N(u) = u-P_N(u)$ for some $P_N(u) \in \bar{P}_N(u)$.  $P_N$ is a ``nonlinear projection'' onto the $N$ ``largest wavelet components'' of a function, and $Q_N$ is the remainder after such a projection.
\begin{lemma}\label{lemma1}
There exists $C=C(p,d) >0$ such that for any $M>0$ and $N\in \N$,
$$\sup_{\scriptsize \begin{array}{c} \|u\|_\lp \leq M,\\ \phantom{sss}Q_N(u) \in \bar{Q}_N(u) \end{array}} \|Q_N(u)\|_\binf \leq C \frac{M}{N^{{1/p}}}\ .$$
\end{lemma}
\noindent
\textbf{Proof:} \quad Fix any $u\in \lprd$ with $\|u\|_\lp \leq M$, fix $\{\l_m\}$ satisfying (\ref{g}) and set $Q_N(u) = \sum_{m=N+1}^\infty a_{\l_m} \psi_{\l_m}$.  Note that, by (\ref{d}) and (\ref{g}),
\begin{equation}\label{h}
\|Q_N(u)\|_\binf \simeq |a_{\l_{N+1}}| \leq \min_{1\leq m \leq N} |a_{\l_m}|\ .
\end{equation}
By (\ref{a1}) of the Appendix we know $\lp \hookrightarrow \dot B^0_{p,p}$, so by (\ref{c}) there exists $C = C(p,d) >0$ such that
\begin{equation}\label{i}
\|a_\l \|_{\ell^p(\LL)} \leq C \|u\|_\lprd
\end{equation}
Now combining (\ref{h}) and (\ref{i}) we see that
$$\left( N \cdot \|Q_N(u)\|^p_\binf \right)^\frac{1}{p} \lesssim \left( \sum_{m=1}^N |a_{\l_m}|^p \right)^\frac{1}{p} \leq
\|a_\l \|_{\ell^p(\LL)}\leq C \|u\|_\lprd \leq CM
$$
and the lemma is proved. \hfill $\Box$
\\\\
We remark here that if one uses the embedding $\dot B^{s_{p,a}}_{a,q} \hookrightarrow \dot B^{s_{p,b}}_{b,b}$ with $b:=\max \{a,q\}$ which is a simple application of Bernstein's inequalities (\ref{aac}) and properties of $\ell^r$ spaces, one easily deduces the similar inequality
\begin{equation}\label{decay2}
\sup_{\scriptsize \begin{array}{c} \|u\|_\besa \leq M,\\ \phantom{sss}Q_N(u) \in \bar{Q}_N(u) \end{array}} \|Q_N(u)\|_\binf \leq C \frac{M}{N^{{1/b}}}\ .
\end{equation}
\section{Main Results}
\noindent
Fix $p,q,r \in \R$ such that $2\leq p < q,r \leq \infty$ and set $s_{p,r}:= d(\frac{1}{r} - \frac{1}{p})$.
Then by Proposition \ref{prop:a3} of the Appendix, there exists $\alpha \in (0,1)$ and $C = C(d,p,q,r) >0$ such that
\begin{equation}\label{j}
\|u\|_\bfin \leq C \|u\|^\alpha_\lprd \|u\|^{1-\alpha}_\binf
\end{equation}
for any $u\in \lprd$.  Note that $\lprd \hookrightarrow \bfin \hookrightarrow \binf$ (see Appendix A).
\\\\
\textbf{Proof of Theorem \ref{thm:a}:}
\\\\
Suppose $\|u_n \|_\lprd \leq \bar{C}$, $\|u_n \tn_\lprd \leq \tilde{C}$, for all $n\in \N$ for some $\bar{C},\tilde{C}>0$, and write ${u_n = \sum_{\l \in \LL} a_{\l,n} \psi_\l}$.  Our first step will be to design an iterative process for describing the sequence $\{u_n\}$ which will allow us to establish (\ref{profiles}) - (\ref{orth4}) of the theorem.
\subsection{Extraction/iteration procedure}
\noindent
\ \\
{({\sc Iterate 0:})
\\\\
Recalling (\ref{k}), if $\lim_{ n\to \infty} \|u_n\|_\binf =0$, then by (\ref{j}) we can write $u_n = r_n^0$ with $\lim_{n\to \infty} \|r_n^0\|_\bfin = 0$, and the statement of the theorem follows.
\\\\
We suppose therefore that $\lim_{n\to\infty}\|u_n\|_\binf \neq 0$, and move to the first iterate:
\\\\
\underline{{\sc Iterate 1:}}
\\\\
By Proposition \ref{prop1}, for each $n\in \N$, there exists $\l_n^1 \in \LL$ such that
$$|a_{\l_n^1}|:=|a_{\l_n^1,n}|  = \max_{\l \in \LL} |a_{\l,n}|\ .$$
By (\ref{d}) and our assumption on $\|u_n\|_\binf$, we may pass to a subsequence such that, for some fixed $i_1 \in \{1, \dots , 2^d -1\}$ and $a_1 \neq 0$,
\begin{equation}\label{l}
\l_n^1 = (i_1, j_n^1,k_n^1) \quad \forall \ n \in \N \qquad \textrm{and} \qquad a_{\l_n^1} \xrightarrow[n\to\infty]{} a_1\ .
\end{equation}
For what follows, let us fix some notation for convenience:
\\\\
If $\tau(x) = ax - b = \frac{x-B}{A}$ for some $a,A \in \R \backslash \{0\}$, $b,B \in \rd$, denote
\begin{equation}\label{m}
(\tau f)(x):=a^{{d/p}} f(\tau(x)) \qquad \textrm{and} \qquad |\tau|:=|\log A| + |B|\ .
\end{equation}
Note that
\begin{equation}\label{n}
\tau_1 (\tau_2 f) = (\tau_2 \circ \tau_1)f
\end{equation}
for any such $\tau_1, \tau_2$.
\\\\
With this notation, we may write
$$u_n = \tau_{1,n} \left(a_{\l_n^1} \psi^{(i_1)} \right) + u_n^1$$
where $\tau_{1,n}(x) := 2^{j_n^1}x - k_n^1$, and $u_n^1$ is just $u_n$ minus one of its wavelet components (i.e., $a_{\l_n^1} \psi_{\l_n^1}$) whose coefficient has the largest possible modulus.  Therefore $\|u_n^1\tn_\binf \leq \|u_n \tn_\binf \leq \tilde{C}$ by (\ref{d}).
\\\\
Suppose $\lim_{n\to\infty}\|u_n^1\|_\binf = 0$.  Then we may write
$$u_n = \tau_{1,n}\left(a_1 \psi^{(i_1)}\right) + r_n^1$$
where $\lim_{n\to\infty}\|r_n^1\|_\bfin = 0$ by (\ref{j}) as follows:
\\\\
Write $r_n^1 = (a_{\l_n^1} - a_1)\psi_{\l_n^1} + u_n^1$, which is just $u_n$ after replacing the coefficient $a_{\l_n^1}$ of $\psi_{\l_n^1}$ by $a_{\l_n^1} - a_1$ which is small for large $n$ by (\ref{l}).  Using (\ref{d}) and (\ref{e}), we therefore see that
$$\|r_n^1 \tn_\binf \leq \|u_n^1 \tn_\binf \xrightarrow[n\to\infty]{} 0$$
and
$$\|r_n^1 \tn_\lprd \leq \|u_n^1 \tn_\lprd \leq \tilde{C}$$
for large $n$.  The statement now follows from (\ref{j}).
\\\\
We have now shown that the statement of the theorem follows when $\lim_{n\to\infty}\|u_n^1\|_\binf = 0$.  We suppose therefore that ${\lim_{n\to\infty}\|u_n^1\|_\binf \neq 0}$, and move to the next iterate:
\\\\
{\sc\underline{Iterate 2:}}
\\\\
By Proposition \ref{prop1} applied to $u_n^1$, for each $n\in \N$, there exists $\l_n^2 \in \LL \backslash \{\l_n^1\}$ such that
$$|a_{\l_n^2}| := |a_{\l_n^2,n}|= \max_{\l \in \LL \backslash \{\l_n^1\}} |a_{\l,n}|\ .$$
By (\ref{d}) and our assumption on $\|u_n^1\|_\binf$, we may pass to a subsequence such that, for some fixed $i_2 \in \{1, \dots , 2^d -1\}$ and $a_2 \neq 0$,
\begin{equation}\label{o}
\l_n^2 = (i_2, j_n^2,k_n^2) \quad \forall \ n \in \N \qquad \textrm{and} \qquad a_{\l_n^2} \xrightarrow[n\to\infty]{} a_2\ .
\end{equation}
Setting $\tau_{2,n}(x):=2^{j_n^2}x - k_n^2$, we can write
$$u_n = \tau_{1,n} \left(a_{\l_n^1} \psi^{(i_1)} \right) + \tau_{2,n} \left(a_{\l_n^2} \psi^{(i_2)} \right) +u_n^2\ .$$
Suppose that $\lim_{n\to\infty}\|u_n^2\|_\binf = 0$.  Let $\tau_n^{(2,1)}:= \tau_{2,n} \circ \tau_{1,n}^{-1}$ and note the following two possibilities (recall (\ref{m})):
\begin{itemize}
\item[(a)] $|\tau_n^{(2,1)}| \to +\infty$ as $n\to\infty$, i.e., ``$\tau_{n,1}$ and $\tau_{n,2}$ are orthogonal'' and
$$u_n = \tau_{1,n} \left(a_1 \psi^{(i_1)} \right) + \tau_{2,n} \left(a_2 \psi^{(i_2)} \right) +r_n^2$$
with $\lim_{n\to\infty}\|r_n^2\|_\bfin = 0$ by (\ref{j}), (\ref{l}) and (\ref{o}) in a way similar to the estimate for $r_n^1$ in {\em Iterate 1}. Note that this means the profiles $a_1 \psi^{(i_1)}$ and  $a_2 \psi^{(i_2)}$ ``live far" (in distance of supports or in scale) from each other after the transformations $\tau_{1,n}$ and $\tau_{2,n}$ respectively, and hence any interactions between the transformed profiles become negligible as $n\to\infty$.
\item[(b)] After possibly passing to a subsequence, $|\tau_n^{(2,1)}| \leq C < \infty$ for all $n$.  We calculate:
$$\tau_n^{(2,1)}(x) = \tau_{2,n}(2^{-j_n^1}(x + k_n^1)) = 2^{j_n^2}(2^{-j_n^1}(x+ k_n^1))-k_n^2 =$$$$= 2^{j_n^2-j_n^1} x - (k_n^2 - 2^{j_n^2-j_n^1} k_n^1) = \frac{x - (2^{j_n^1-j_n^2}k_n^2 - k_n^1)}{2^{j_n^1-j_n^2}}\ ,$$
so (\ref{m}) gives
$$|\log 2^{j_n^1-j_n^2}| + |2^{j_n^1-j_n^2}k_n^2 - k_n^1| \leq C\ .$$
Note therefore that $j_n^2 - j_n^1 \equiv j^{(2,1)}$ for infinitely many $n$ for some $j^{(2,1)} \in \Z$, and then $k_n^2 - 2^{j^{(2,1)}} k_n^1 \equiv b^{(2,1)}$ for infinitely many of those $n$ for some $b^{(2,1)} \in \rd$, due to the fact that we are working on lattices.  We may therefore pass to a subsequence so that
$$\tau_n^{(2,1)} \equiv \tau^{(2,1)}= \textrm{const} \quad \textrm{for all} \ n\in \N, $$
with
$$\tau^{(2,1)}(x) = 2^{j^{(2,1)}}x - b^{(2,1)}\ .$$
Recalling (\ref{n}), we may now write
$$
\begin{array}{rcl}
u_n & = & \tau_{1,n} \left(a_{\l_n^1} \psi^{(i_1)} \right) + \tau_{2,n} \left(a_{\l_n^2} \psi^{(i_2)} \right) +u_n^2\\\\
 & = & \tau_{1,n} \left( a_{\l_n^1} \psi^{(i_1)} + \tau_{1,n}^{-1} (\tau_{2,n} (a_{\l_n^2} \psi^{(i_2)})) \right) +u_n^2 \\\\
 & = & \tau_{1,n} \left( a_{\l_n^1} \psi^{(i_1)} + (\tau_{2,n} \circ \tau_{1,n}^{-1}) (a_{\l_n^2} \psi^{(i_2)}) \right) +u_n^2 \\\\
 & = & \tau_{1,n} \left( a_{\l_n^1} \psi^{(i_1)} + a_{\l_n^2} \tau^{(2,1)}  \psi^{(i_2)} \right) +u_n^2 \\\\
 & =: & \tau_{1,n} \left( a_1 \psi^{(i_1)} + a_2 \tau^{(2,1)}  \psi^{(i_2)} \right) +r_n^2\ .
\end{array}
$$
Comparing the last two lines, we can again show that ${\lim_{n\to\infty} \|r_n^2\|_\bfin = 0}$ in a similar way to {\em Iterate 1} using (\ref{j}), (\ref{l}) and (\ref{o}).
\end{itemize}
This concludes {\em Iterate 2}, and the statement of the theorem follows if ${\lim_{n\to\infty}\|u_n^2\|_\bfin = 0}$.  If not, we would move on to the next iterate.
\\\\
Let us however pause here (before giving the full induction argument for further iterates) to give some geometric meaning to this dichotomy.  For $m=1,2$, note that $|\d_n^m| = (2^{-j_n^m})^d$, and $\d_n^m$ (which corresponds to the location of $\psi_{\l_n^m}$) has $2^{-j_n^m}k_n^m$ as its ``center" (more precisely, the location of the ``corner'' located closest to the origin).
\\\\
Possibility (a) means that either
\begin{itemize}
\item[(i)] the sizes of the cubes $\d_n^1$ and $\d_n^2$ are incomparable as $n\to\infty$; that is,
$$\left|\log \frac{2^{-j_n^2}}{2^{-j_n^1}}\right| = | \log 2^{j_n^1 - j_n^2}| \xrightarrow[n\to\infty]{} \infty\ ,$$
or
\item[(ii)] they have comparable size $2^{-j_n^2} \simeq 2^{-j_n^1}$ for all $n$, and the distance between the cubes {\em relative to their size} becomes infinite as $n\to\infty$.  That is, after dilating space by, say, $2^{j_n^1}$, bringing their sizes to $O(1)$, their new centers ($k_n^1$ and $2^{j_n^1 - j_n^2}k_n^2$) separate:
    $$\left| \frac{2^{-j_n^2}k_n^2 - 2^{-j_n^1}k_n^1}{2^{-j_n^1}} \right| = |2^{j_n^1 - j_n^2}k_n^2 - k_n^1| \xrightarrow[n\to\infty]{} \infty\ .$$
\end{itemize}
In either case, we see that $\psi_{\l_n^1}$ and $\psi_{\l_n^2}$ have negligible interaction for large $n$.
\\\\
Possibility (b) means neither (i) nor (ii) happen, so the cubes have comparable size and {\em relative} (to their size) distance.  Note that $\tau_{m,n}\d_n^m = [0,1)^d$, and we are saying that for a subsequence, after rescaling to $O(1)$ by, say, $\tau_{n,1}$, the relative positions of the cubes are
constant:\footnote{Note that this was stated slightly inaccurately in \cite{jaffard} -- we hope that this explanation clarifies the matter.}}
$$\tau^{(2,1)}(\tau_{1,n} \d_n^2) = \tau_{1,n}\d_n^1\ .$$
This gives $\tau^{(2,1)} \tau_{1,n} = \tau_{2,n} \iff \tau^{(2,1)} = \tau_{2,n}\tau_{1,n}^{-1}$, as we had previously designated due to {\em algebraic} considerations.
\\\\
We will now continue to inductively define all possible iterates.  From now on, when we pass to a subsequence at the $N^{\textrm{th}}$ iterate, we will leave the first $N$ sequence elements unchanged, obtaining a ``diagonal'' subsequence which will work for all iterates.
\\\\
Suppose, after $N$ iterations, $\lim_{n\to\infty} \|u_n^{N'}\|_\binf \neq 0$ for all $N' < N$, and the conclusion of the $N^{\textrm{th}}$ iterate is as follows (and holds for all previous iterates):
\\\\
{\sc \underline{Iterate N:}}
\begin{equation}\label{p}
u_n = \sum_{l=1}^{L_N} \tau_{m_1(l),n} \left( \sum_{\mu = 1}^{M_N(l)} a_{\l_n^{m_\mu(l)}} \tau^{(m_\mu(l),m_1(l))} \psi^{(i_{m_\mu(l)})}  \right) + u_n^N
\end{equation}
where $\sum_{l=1}^{L_N} M_N(l) = N$, $m_\mu(l)$ are never repeated (i.e., $m_\mu(l) = m_{\mu'}(l') \Longrightarrow {\mu = \mu', l=l'}$), $i_{m_\mu(l)} \in \{1,\dots, 2^d -1\}$,
$$\tau_{m_\mu(l),n}(x) := 2^{j_n^{m_\mu(l)}}x - k_n^{m_\mu(l)}$$
and
$$\tau^{(m_\mu(l),m_1(l))}(x): = 2^{j^{(m_\mu(l),m_1(l))}}x - b^{(m_\mu(l),m_1(l))}$$
for some $j_n^{m_\mu(l)},j^{(m_\mu(l),m_1(l))} \in \Z$, $k_n^{m_\mu(l)}\in \Z^d$ and $b^{(m_\mu(l),m_1(l))} \in \rd$ such that the following hold:
\begin{enumerate}
\item For each $\mu$ and $l$,
$$\l_n^{m_\mu(l)} = (i_{m_\mu(l)}, j_n^{m_\mu(l)},k_n^{m_\mu(l)}) \quad \forall \ n \in \N$$
\item $j_n^{m_\mu(l)} - j_n^{m_1(l)} \equiv j^{(m_\mu(l),m_1(l))}= \textrm{const}$ for all $n\in \N$,
\item $k_n^{m_\mu(l)} - 2^{j^{(m_\mu(l),m_1(l))}}k_n^{m_1(l)} \equiv b^{(m_\mu(l),m_1(l))}= \textrm{const}$ for all $n\in \N$,
\item $\tau_{m_\mu(l),n} \tau_{m_1(l),n}^{-1} \equiv \tau^{(m_\mu(l),m_1(l))} = \textrm{const}$ for all $n\in \N$,\\\\ (4. is a consequence of 2. and 3.)
\item for any $l\neq l'$,
\begin{equation}\label{six}
|\tau_{m_1(l),n}\tau_{m_1(l'),n}^{-1}|  \xrightarrow[n\to\infty]{} +\infty\ ,
\end{equation}
\item $\displaystyle{|a_{\l_n^m}|:= |a_{\l_n^m,n}| = \max_{\l \in \LL \backslash \{\l_n^{m'}\}_{m'=1}^{m-1}} |a_{\l,n}|}$ for any $m \in \N$, and
\item there exists $a_{m_\mu(l)} \in \R \backslash \{0\}$ for each $\mu$ and $l$ such that
\begin{equation}\label{q}
a_{\l_n^{m_\mu(l)}} \xrightarrow[n\to\infty]{} a_{m_\mu(l)}\ .
\end{equation}
\end{enumerate}
Note that in {\em Iterate 2}, case (a) we had the above with $L_2 =2$, $M_2(1) = M_2(2) =1$, $m_1(1)=1$ and $m_1(2) =2$, and {\em Iterate 2} case (b) was $L_2 =1$, $M_2(1) =2$, $m_1(1) =1$, $m_2(1) =2$.
\\\\
If $\lim_{n\to\infty}\|u_n^N\|_\binf =0$, then the statement of the theorem follows as before by writing
\begin{equation}\label{r}
u_n = \sum_{l=1}^{L_N} \tau_{m_1(l),n} \left( \sum_{\mu = 1}^{M_N(l)} a_{{m_\mu(l)}} \tau^{(m_\mu(l),m_1(l))} \psi^{(i_{m_\mu(l)})}  \right) + r_n^N\ ,
\end{equation}
comparing (\ref{r}) with (\ref{p}) and estimating $r_n^N$ as in {\em Iterate 1} and {\em Iterate 2} using (\ref{j}) and (\ref{q}).
\\\\
Therefore we assume again that $\lim_{n\to\infty}\|u_n^N\|_\binf \neq 0$, and move to the next iterate:
\\\\
{\sc \underline{Iterate N+1:}}
\\\\
We pass to a subsequence so that
$$a_{\l_n^{N+1}} \xto a_{N+1} \neq 0$$
where $|a_{\l_n^{N+1}}|:=|a_{\l_n^{N+1},n}| = \|u_n^N \tn_\binf$ for each $n$ and
$$\l^{N+1}_n = (i_{N+1},j^{N+1}_n,k_n^{N+1})\in \LL$$
for all $n$ for some fixed $i_{N+1} \in \{1,\dots,2^d -1\}$.
\\\\
Let $\tau_{N+1,n}(x):= 2^{j^{N+1}_n}x-k_n^{N+1}$ and define $\tau_n^{(N+1,m_1(l))}:=\tau_{{N+1},n}\tau_{m_1(l),n}^{-1}$ for $l = 1, \dots, L_N$.  We then have the following two possibilities:
\begin{itemize}
\item[(a)]  $|\tau_n^{(N+1,m_1(l))}| \xto +\infty$ for all $l\in \{1,\dots,L_N\}$.  In this case, we set
$$L_{N+1} = L_N + 1, \quad M_{N+1}(L_{N+1}) = 1 \quad \textrm{and} \quad m_1(L_{N+1}) = N+1\ $$
which gives (\ref{p}) with $N$ replaced by $N+1$ with all the required properties.  That is, we start a new function in slot $L_{N}+1$ so that $u_n$ is a composed of rescalings of $L_N + 1$ functions (as opposed to $L_N$ functions as in {\em Iterate N}), the new function being $a_{\l_n^{N+1}} \psi^{(i_{N+1})}$.  (These functions become constant in $n$ after writing this in the form of (\ref{r}), the new function being $a_{N+1} \psi^{(i_{N+1})}$.)
\item[(b)]  There is some $\bar{l} \in \{1,\dots,L_N\}$ and $C>0$ such that, for some subsequence,
$$|\tau_n^{(N+1,m_1(\bar{l}))}| \leq C \quad \textrm{for all}\ n\in \N\ .$$
In this case, we pass to a subsequence so that
$$\tau_n^{(N+1,m_1(\bar{l}))} \equiv \tau^{(N+1,m_1(\bar{l}))} = \textrm{const} \quad \textrm{for all}\ n\in\N$$
with $\tau^{(N+1,m_1(\bar{l}))}$ defined in a way analogous to 2., 3. and 4. of ${\textit{Iterate N}}$ with $m_\mu(l)$ replaced by $N+1$ and $l$ replaced by $\bar{l}$ everywhere else.  We then set
$$L_{N+1} = L_N, \quad M_{N+1}(\bar{l}) = M_N(\bar{l}) + 1 \quad \textrm{and} \quad m_{M_{N+1}(\bar{l})} = N+1\ .$$
That is, we still have re-scalings of $L_N$ functions as in {\em Iterate N}, and we just add the new component, $a_{\l_n^{N+1}} \psi_{\l_n^{N+1}}$, to the $\bar{l}^{\textrm{th}}$ one so that the $\bar{l}^{\textrm{th}}$ term in the sum in (\ref{p}) becomes
$$\tau_{m_1(\bar{l}),n} \left( \left( \sum_{\mu = 1}^{M_N(\bar{l})} a_{\l_n^{m_\mu(\bar{l})}} \tau^{(m_\mu(\bar{l}),m_1(\bar{l}))} \psi^{(i_{m_\mu(\bar{l})})} \right)+ a_{\l_n^{N+1}} \tau^{(N+1,m_1(\bar{l}))} \psi^{(i_{N+1})} \right)\ ,$$
and all others remain the same.
\end{itemize}
In either case, if $\lim_{n\to\infty}\|u_n^{N+1}\|_\binf = 0$, the statement of the theorem follows as in {\em Iterate N}.  If $\lim_{n\to\infty}\|u_n^{N+1}\|_\binf \neq 0$ then we move to {\em Iterate N+2} and so on, so we have now inductively defined all possible iterates.
\subsection{Proof of properties (\ref{profiles}) - (\ref{orth2}) and (\ref{orth4})}
If $\lim_{n\to\infty}\|u_n^{\bar{N}}\|_\binf = 0$ for some $\bar{N} \in \N$ then we stop the iterations at ${\mathit{Iterate \bar{N}}}$ and the main statement of the theorem follows (with $r_n^L \equiv r_n^{\bar{N}}$, etc., for large $L$).
\\\\
Therefore we suppose now that $\lim_{n\to\infty}\|u_n^N\|_\binf \neq 0$ for all $N\in \N$, and we pass to a diagonal subsequence so that the statements in {\em Iterate N} hold for all $N\in \N$ (without changing subsequences).
\\\\
Let $L_\infty := \lim_{N\to \infty}L_N$ and $M_\infty(l):= \lim_{N\to \infty}M_N(l)$ for any $l\in [1,L_\infty] \cap \N$.  Then for any $L\in [1,L_\infty] \cap \N$ we write
$$u_n = \sum_{l=1}^{L} \tau_{m_1(l),n} \left( \sum_{\mu = 1}^{M_\infty(l)} a_{{m_\mu(l)}} \tau^{(m_\mu(l),m_1(l))} \psi^{(i_{m_\mu(l)})}  \right) + r_n^L
$$
\begin{equation}\label{s}
\qquad \qquad  = \sum_{l=1}^{L}  \left( \sum_{\mu = 1}^{M_\infty(l)} a_{{m_\mu(l)}} \tau_{m_\mu(l),n} \psi^{(i_{m_\mu(l)})}  \right) + r_n^L  \ .
\end{equation}
If $M_\infty(l) < \infty$ for some $l$, the term in the parentheses in (\ref{s}) is a well-defined element of $\lprd$, since it is a wavelet expansion with only finitely many components.  Hence we have also
\begin{equation}\label{aa}
\phi_l := \sum_{\mu = 1}^{M_\infty(l)} a_{{m_\mu(l)}} \tau^{(m_\mu(l),m_1(l))} \psi^{(i_{m_\mu(l)})} \in \lprd
\end{equation}
by the invariance of the usual $\lprd$ norm under $ \tau_{m_1(l),n}$.
\\\\
If $M_\infty(l) = +\infty$ then (\ref{e}) and Fatou's lemma (along with the boundedness of $\{u_n\}$) imply $\phi_l \in \lprd$ as well, a fact we leave to the Appendix.  Indeed, apply Lemma \ref{lemma:a4} of the Appendix with $a_n^\mu = a_{\l_n^{m_\mu(l)}}$, $\l_n^\mu = \l_n^{m_\mu(l)}$ and $C=\bar{C}$.
\\\\
Property (\ref{orth1}) is clear from our construction, so we must now show that
$$\lim_{L \to L_\infty} \left( \varlimsup_{n\to\infty} \|r_n^L\|_\bfin \right) =0$$
which will establish (\ref{orth2}), for which we turn to Lemma \ref{lemma1}.
\\\\
Let us first assume that $L_\infty = +\infty$ (the case $L_\infty < +\infty$ is even easier).
\\\\
Note first by standard embeddings (see the Appendix) that we have
$$\sum_{m=1}^\infty |a_{\l_n^m}|^p = \|u_n\tn_{\dot B^0_{p,p}}^p \lesssim \|u_n\|_\lprd^p \leq \bar{C}^p \quad \textrm{for all}\ n\in\N\ ,$$
so by Fatou's lemma we have
$$\sum_{m=1}^\infty |a_m|^p = \sum_{m=1}^\infty \varliminf_{n\to\infty}|a_{\l_n^m}|^p \leq \varliminf_{n\to\infty}\sum_{m=1}^\infty |a_{\l_n^m}|^p \lesssim \bar{C}^p < \infty\ ,$$
hence
\begin{equation}\label{t}
\lim_{m\to\infty}|a_m| =0\ .
\end{equation}
We will need this fact in what follows.
\\\\
Fix any $L\in \N$.  For any $N$ large enough that $L_N \geq L$ (recall we are assuming $L_\infty = +\infty$), by (\ref{p}) we can write
$$u_n = \sum_{l=1}^{L_N}  \left( \sum_{\mu = 1}^{M_N(l)} a_{\l_n^{m_\mu(l)}} \tau_{m_\mu(l),n} \psi^{(i_{m_\mu(l)})}  \right) + u_n^N$$
$$= \sum_{l=1}^{L}  \left( \sum_{\mu = 1}^{M_\infty(l)} a_{m_\mu(l)} \tau_{m_\mu(l),n} \psi^{(i_{m_\mu(l)})}  \right) + r_n^L$$
which gives
\begin{equation}\label{u}
\begin{array}{rcl}
r_n^L & = & \displaystyle{\sum_{l=1}^{L}  \left( \sum_{\mu = 1}^{M_N(l)} (a_{\l_n^{m_\mu(l)}} - a_{m_\mu(l)}) \tau_{m_\mu(l),n} \psi^{(i_{m_\mu(l)})}  \right)}\\\\
&  & \quad - \ \displaystyle{ \sum_{l=1}^{L} \left( \sum_{\mu = M_N(l)+1}^{M_\infty(l)} a_{m_\mu(l)} \tau_{m_\mu(l),n} \psi^{(i_{m_\mu(l)})}  \right)      }\\\\
&  & \quad + \ \displaystyle{ \sum_{l=L+1}^{L_N} \left( \sum_{\mu = 1}^{M_N(l)} a_{\l_n^{m_\mu(l)}} \tau_{m_\mu(l),n} \psi^{(i_{m_\mu(l)})}  \right) \ +\  u_n^N     }\\\\
& =: & I - II + III + IV
\end{array}
\end{equation}
(where if $L=L_N$ or $M_N(l) = M_\infty(l)$, we replace the appropriate sum by zero).  We'll estimate, using the embedding $\lprd \hookrightarrow \bfin$,
\begin{equation}\label{v}
\|r_n^L\|_\bfin \lesssim \|I\|_\lprd + \|II\|_\lprd + \|III\|_\bfin + \|IV\|_\bfin\ .
\end{equation}
We claim, in fact, that this leads to an estimate of the form
\begin{equation}\label{w}
\|r_n^L\|_\bfin \leq \e_{N,L}(n) + \e_L^1(N) + \e_L^2(n) + \frac{C}{N^\frac{1}{p}}\ ,
\end{equation}
where
\begin{equation}\label{x}
\begin{array}{c}
\displaystyle{\lim_{n\to\infty} \e_{N,L}(n) = 0 \ \ \forall\ N,L \in \N\ , \quad  \lim_{N\to\infty} \e_L^1(N) =0 \ \ \forall\ L \in \N\ ,}\\\\
\displaystyle{\textrm{and}\qquad \lim_{L\to\infty}\e_L^2=0 \quad \textrm{where} \quad \e_L^2:=\lim_{n\to\infty} \e_L^2(n)\ .}
\end{array}
\end{equation}
In this case, the desired estimate on $r_n^L$ is proved as follows:  For any $\bar{\e}>0$, there exists $\bar{L}\in \N$ such that $\e^2_{L_0} < {\bar{\e}/3}$ for any $L_0 \geq \bar{L}$.  For any such $L_0$, there exists $N_0 = N_0(L_0)\in \N$ large enough that $L_{N_0} \geq L_0$ and ${\max \{\e_{L_0}^1(N_0), {{C N_0^{-{1/p}}}} \} < {\bar{\e}/3}}$.  Hence by (\ref{w})-(\ref{x}), for any $\bar{\e}>0$ there exists $\bar{L}\in \N$ such that, for all $L_0 \geq \bar{L}$,
$$\varlimsup_{n\to\infty} \|r_n^{L_0}\|_\bfin \leq \lim_{n\to\infty} \left( \e_{N_0,L_0}(n)\ +\ \e_{L_0}^1(N_0)\ +\ \e_{L_0}^2(n)\ +\ \frac{C}{N_0^\frac{1}{p}}\right) =$$
$$\qquad \qquad \qquad \qquad
\qquad = \qquad \  0\quad  +\  \e_{L_0}^1(N_0) \ +\  \e_{L_0}^2\  +\  \frac{C}{N_0^\frac{1}{p}}\qquad   <\qquad  \bar{\e}\ ,
$$
which is exactly (\ref{orth2}).
\\\\
Let us pause here to address (\ref{orth4}) as well:  It is clear that $III+IV$ is just $u_n$ with some of its wavelet components removed.  Therefore, using (\ref{e}) and (\ref{u}) - (\ref{w}), we can estimate
$$\|r_n^L\tn_\lprd \leq \e_{N,L}(n) + \e_L^1(N) + \|u_n\tn_\lprd
$$
and hence for any $\e>0$, by taking first $N$ large and then $n$ large, we see by (\ref{x}) that for any $L\in \N$ we have
$$\left| \|r_n^L\tn_\lprd - \|u_n\tn_\lprd  \right| < \e $$
for $n$ sufficiently large which establishes (\ref{orth4}).
\\\\
Turning back to the proof of (\ref{orth2}), we'll now estimate individually each term of (\ref{v}).  Note first the following important change of variable formula:
\\\\
For any $l$ and $M$ and any numbers $a_\mu (n)$ possibly depending on $n$, we may write
\begin{equation}\label{y}
\sum_{\mu=1}^M |a_\mu(n)|^2 2^{\frac{2d}{p} j_n^{m_\mu(l)}} \chi_{\d_n^{m_\mu(l)}}(x) = 2^{\frac{2d}{p}j_n^{m_1(l)}} \tilde{\chi}_l (2^{j_n^{m_1(l)}}x - k_n^{m_1(l)} )
\end{equation}
where
\begin{equation}\label{z}
\tilde{\chi}_l(y) = \sum_{\mu=1}^M |a_\mu(n)|^2 2^{\frac{2d}{p}j^{(m_\mu(l),m_\mu(1))}} \chi_{\d^{(m_\mu(l),m_\mu(1))}}(y)
\end{equation}
and
$$\d^{(m_\mu(l),m_\mu(1))} = (\tau^{(m_\mu(l),m_\mu(1))})^{-1} ([0,1)^d)$$
(see the proof of Lemma \ref{lemma:a4}).  Recall that $\chi_\d$ is the characteristic function of the cube $\d$, and note that $\tilde \chi_l$ is bounded with compact support for each $l$, and of course is independent of $n$.
\pagebreak
\\\\
\underline{{\sc Estimate of I:}}
\\\\
Using (\ref{e}) and (\ref{y}), for any $l\in \{1,\dots,L\}$ the change of variable $y=2^{j_n^{m_1(l)}}x - k_n^{m_1(l)}$ gives
$$\Bigg\| \sum_{\mu=1}^{M_N(l)} (a_{\l_n^{m_\mu(l)}} - a_{m_\mu(l)}) \tau_{m_\mu(l),n} \psi^{(i_{m_\mu(l)})}
\tilde{\Bigg\|}^p_\lprd =$$
$$=
\int \left\{ \sum_{\mu=1}^{M_N(l)}|a_{\l_n^{m_\mu(l)}} - a_{m_\mu(l)}|^2 2^{\frac{2d}{p}j^{(m_\mu(l),m_\mu(1))}} \chi_{\d^{(m_\mu(l),m_\mu(1))}}(y) \right\}^\frac{p}{2}dy
$$
$$
\leq 2 \int \left\{ \sum_{\mu=1}^{M_N(l)}|a_{m_\mu(l)}|^2 2^{\frac{2d}{p}j^{(m_\mu(l),m_\mu(1))}} \chi_{\d^{(m_\mu(l),m_\mu(1))}}(y) \right\}^\frac{p}{2}dy \leq \tilde{C}^p < \infty
$$
for all $n$ sufficiently large, depending on $N$, and the last inequality is clear from the proof of Lemma \ref{lemma:a4}.  Lebesgue's Dominated Convergence Theorem now gives that the integral on the left tends to zero as $n\to \infty$ for fixed $N$, hence, after summing from $1$ to $L$ ($L<\infty$), we have
$$\|I\|_\lprd \leq \e_{N,L}(n) \qquad \textrm{with} \qquad \lim_{n\to\infty} \e_{N,L}(n) = 0$$
as desired.
\\\\
\underline{{\sc Estimate of II:}}
\\\\
Similarly, for $l \in \{1,\dots,L\}$, we have
$$\Bigg\| \sum_{\mu=M_N(l)+1}^{M_\infty(l)} a_{m_\mu(l)} \tau_{m_\mu(l),n} \psi^{(i_{m_\mu(l)})}
\tilde{\Bigg\|}^p_\lprd =$$
$$=
\int \left\{ \sum_{\mu=M_N(l)+1}^{M_\infty(l)} |a_{m_\mu(l)}|^2 2^{\frac{2d}{p}j^{(m_\mu(l),m_\mu(1))}} \chi_{\d^{(m_\mu(l),m_\mu(1))}}(y) \right\}^\frac{p}{2}dy
$$
$$
\leq  \int \left\{ \sum_{\mu=1}^{M_\infty(l)}|a_{m_\mu(l)}|^2 2^{\frac{2d}{p}j^{(m_\mu(l),m_\mu(1))}} \chi_{\d^{(m_\mu(l),m_\mu(1))}}(y) \right\}^\frac{p}{2}dy \leq \tilde{C}^p < \infty
$$
uniformly in $N$ by Lemma \ref{lemma:a4}.  Now, again by Lebesgue's Dominated Convergence Theorem, the left-hand side tends to zero as $N\to \infty$.  Summing from $1$ to $L$, we have
$$\|II\|_\lprd \leq \e_L(N) \qquad \textrm{where} \qquad \lim_{N\to\infty} \e_L(N) =0$$
as desired.
\\\\
\underline{{\sc Estimate of III:}}
\\\\
Since $III$ is just $u_n$ after replacing some of its wavelet coefficient by zero, (\ref{e}) gives
\begin{equation}\label{three}
\|III\tn_\lprd \leq \|u_n \tn_\lprd \leq \tilde{C}\ .
\end{equation}
Note that, by our construction, $m_1(l)$ is increasing with $l$ and $m_\mu(l)$ is increasing with $\mu$ for each fixed $l$.  Therefore, since also $|a_{\l_n^m}|$ is monotonically decreasing in $m$ for each fixed $n$, we see that
$$\|III\tn_\binf = |a_{\l_n^{m_1(L+1)}}|$$
which is independent of $N$.  Since $m_1(L+1) \to +\infty$ as $L\to \infty$, (\ref{t}) gives
\begin{equation}\label{four}
\|III\tn_\binf \leq \e_L(n)\ ,
\end{equation}
where
\begin{equation}\label{five}
\lim_{L\to\infty} \left(\lim_{n\to\infty} \e_L(n) \right) = \lim_{L\to\infty} |a_{m_1(L+1)}| = 0\ .
\end{equation}
The desired estimate now follows from (\ref{three}) - (\ref{five}) and the interpolation inequality (\ref{j}).
\\\\
\underline{{\sc Estimate of IV:}}
\\\\
Similarly, $\|u_n^N\tn_\lprd \leq \|u_n \tn_\lprd \leq \tilde{C}$, and since $u_n^N \in \bar{Q}_N(u_n)$, Lemma \ref{lemma1} gives $\|u_n^N\|_\binf \leq CN^{-{{1/p}}}$ for some $C>0$, independently of $n$.  The desired estimate now follows from (\ref{j}).
\\\\
Estimate (\ref{w}) is now proved, and defining $\phi_l$ as in (\ref{aa}), we have
$$u_n = \sum_{l=1}^L \tau_{m_1(l),n} \phi_l + r_n^L, \qquad \textrm{with} \qquad \lim_{L\to\infty} \left( \varlimsup_{n\to\infty} \|r_n^L\|_\bfin \right) =0$$
thus establishing (\ref{profiles}) - (\ref{orth2}) of the theorem in the case $L_\infty = \infty$.  The case $L_\infty < \infty$ is even easier; starting with $L=L_\infty$ from the beginning, the term $III$ completely disappears from (\ref{u}) so that $\lim_{n\to\infty} \|r_n^{L_\infty}\|_\bfin = 0$ by the considerations above.
\subsection{Proof of (\ref{orth3})}
We now turn to the proof of (\ref{orth3}):
\\\\
For convenience, let's define
$$\l_n^0 = (0,0,0)\ , \quad a_{\l_n^0} = a_0 = \psi^{(0)} = j^{(0,m_1(l))} = b^{(0,m_1(l))} =0\ ,$$
and replace both $L_\infty$ and $M_\infty(l)$ by $\infty$ by making the convention that in case either is finite, we set $m_\mu(l) = 0$ whenever $l> L_\infty$ or $\mu > M_\infty(l)$.
\\\\
With this convention, we can write
\begin{equation}\label{ab}
u_n = \sum_{l=1}^{\infty}  \sum_{\mu = 1}^{\infty} a_{\l_n^{m_\mu(l)}} \tau_{m_\mu(l),n} \psi^{(i_{m_\mu(l)})}  \ ,
\end{equation}
which is just a particular ordering (depending on $n$), of the wavelet components of $u_n$.  Therefore (\ref{e}) gives
\begin{equation}\label{ac}
\int \left\{ \sum_{l=1}^{\infty}  \sum_{\mu = 1}^{\infty} |a_{\l_n^{m_\mu(l)}}|^2 2^{\frac{2d}{p}j_n^{m_\mu(l)}} \chi_{\d_n^{m_\mu(l)}} (x) \right\}^\frac{p}{2} dx = \|u_n \tn_\lprd^p \leq \tilde{C}^p
\end{equation}
uniformly in $n$.  A change of variables as in (\ref{y}) and an application of Lemma \ref{lemma:a4} shows that for any $l\in \N$ (recall (\ref{aa})),
\begin{equation}\label{za}
\begin{array}{c}
\displaystyle{\|\tau_{m_1(l),n}\phi_l \tn_p^p = \left\{   \sum_{\mu = 1}^{\infty} |a_{m_\mu(l)}|^2 2^{\frac{2d}{p}j_n^{m_\mu(l)}} \chi_{\d_n^{m_\mu(l)}} (x) \right\}^\frac{p}{2} dx  = \qquad } \\\\
\phantom{f}  \qquad \displaystyle{=\left\{   \sum_{\mu = 1}^{\infty} |a_{m_\mu(l)}|^2 2^{\frac{2d}{p}j^{(m_\mu(l),m_\mu(1))}} \chi_{\d^{(m_\mu(l),m_\mu(1))}}(y) \right\}^\frac{p}{2}dy \leq   \tilde{C}^p}
\end{array}
\end{equation}
for any $n\in\N$ -- that is, within each profile one can replace the wavelet coefficients by their pointwise limits.  What we want to show now is that
$$\sum_{l=1}^\infty \|\tau_{m_1(l),n}\phi_l \tn_p^p \leq \liminf_{n\to\infty} \|u_n\tn_p^p, $$
which is (\ref{orth3}), and at this point we have established that all quantities are well-defined.
\\\\
To simplify notation, let's denote
$$\eta^{\mu,l}_n(x):= 2^{\frac{2d}{p}j_n^{m_\mu(l)}} \chi_{\d_n^{m_\mu(l)}} (x) \quad \textrm{and} \quad \tilde{\eta}^{\mu,l}(y):= 2^{\frac{2d}{p}j^{(m_\mu(l),m_\mu(1))}} \chi_{\d^{(m_\mu(l),m_\mu(1))}}(y)\ .$$
Truncating (\ref{ac}) and (\ref{za}) to finite sums, Fatou's lemma applied as in Lemma \ref{lemma:a4} yields
$$\sum_{l=1}^L \int \left\{   \sum_{\mu = 1}^{M} |a_{m_\mu(l)}|^2 \tilde{\eta}^{\mu,l}(y) \right\}^\frac{p}{2}dy
\leq \liminf_{n\to\infty} \sum_{l=1}^L \int \left\{   \sum_{\mu = 1}^{M} |a_{\l_n^{m_\mu(l)}}|^2 \tilde{\eta}^{\mu,l}(y) \right\}^\frac{p}{2}dy .
$$
We claim moreover that
$$\sum_{l=1}^L \int \left\{   \sum_{\mu = 1}^{M} |a_{\l_n^{m_\mu(l)}}|^2 \eta^{\mu,l}_n(x) \right\}^\frac{p}{2}dx
= \int \left\{  \sum_{l=1}^L \sum_{\mu = 1}^{M} |a_{\l_n^{m_\mu(l)}}|^2 \eta^{\mu,l}_n(x) \right\}^\frac{p}{2}dx  +\e_1(n)\ ,
$$
where $\lim_{n\to\infty}\e_1(n) =0$.  (This is simply the quantity on the right of the previous line after a change of variable.)  Let us postpone this claim for the moment.
\\\\
Using (\ref{ac}) to bound the last term above by $\|u_n\tn_p^p + \e_1(n)$ and inserting this into the previous inequality we have
$$\sum_{l=1}^L \int \left\{   \sum_{\mu = 1}^{M} |a_{m_\mu(l)}|^2 \tilde{\eta}^{\mu,l}(y) \right\}^\frac{p}{2}dy \leq \liminf_{n\to\infty}\|u_n\tn_p^p\ .$$
Now that both sides are independent of n, using Fatou's lemma to let $M\to\infty$ and then letting $L\to\infty$ and using identity (\ref{za}) yields (\ref{orth3}) as desired.
\\\\
We are reduced therefore to showing that
\begin{equation}\label{ag}
\begin{array}{c}
    \displaystyle{ \lim_{n\to\infty} I_{L,M}(n) := }\\\\
 \displaystyle{ \lim_{n\to\infty}\Bigg| \int \left\{ \sum_{l=1}^{L}  \sum_{\mu = 1}^{M} |a_{\l_n^{m_\mu(l)}}|^2 2^{\frac{2d}{p}j_n^{m_\mu(l)}} \chi_{\d_n^{m_\mu(l)}} (x) \right\}^\frac{p}{2} dx } \qquad  \qquad  \qquad \\\\
 \qquad  \qquad  \qquad   \displaystyle{ -\ \sum_{l=1}^{L} \int \left\{  \sum_{\mu = 1}^{M} |a_{\l_n^{m_\mu(l)}}|^2 2^{\frac{2d}{p}j_n^{m_\mu(l)}} \chi_{\d_n^{m_\mu(l)}} (x) \right\}^\frac{p}{2} dx} \Bigg| =0
\end{array}
\end{equation}
Note first that if $p=2$, (\ref{ag}) is trivial so we may assume that $p>2$.
\\\\
We'll use the elementary inequality (see, e.g., equation (1.10) in \cite{gerard}) that for any $L\in \N$ there exists a constant $C_L=C_L(p) >0$ such that for any $\{A_l\} \subset [0,\infty)$,
\begin{equation}\label{aj}
\left| \left( \sum_{l=1}^L A_l \right)^\frac{p}{2} - \sum_{l=1}^L (A_l)^\frac{p}{2} \right| \leq \ \ C_L \sum_{l\neq l'} A_l (A_{l'})^{\frac{p}{2}-1}\ .
\end{equation}
We can therefore estimate the function in (\ref{ag}) using (\ref{aj}) and the notation of (\ref{y}) and (\ref{z}) as
$$I_{L,M}(n) \leq C_L \sum_{l\neq l'} I^{l,l'}(n)$$
with
$$I^{l,l'}(n):= \int \left(2^{dj_n^{m_1(l)}}\right)^{(\frac{2}{p})} \tilde \chi_l (2^{j_n^{m_1(l)}}x - k_n^{m_1(l)})
\left(2^{dj_n^{m_1(l')}}\right)^{(1-\frac{2}{p})} \tilde \chi_{l'}^{(\frac{p}{2}-1)} (2^{j_n^{m_1(l')}}x - k_n^{m_1(l')})\ dx\ .
$$
(We have absorbed the coefficients $|a_{\l_n^{m_\mu(l)}}|^2$ into the constant since they converge and hence are uniformly bounded.) Fix $l$ and $l'$, and set the notation
\begin{equation}\label{zzz}
x_n^{l_1,l_2}(w):=2^{j_n^{m_1(l_2)}-j_n^{m_1(l_1)}}w + \left[2^{j_n^{m_1(l_2)}-j_n^{m_1(l_1)}}k_n^{m_1(l_1)}- k_n^{m_1(l_2)}\right]
\end{equation}
for $w\in \rd$ and any $l_1,l_2\in \N$.  Then by appropriate changes of variables, we see
\begin{equation}\label{ak}
\begin{array}{rcl}
I^{l,l'}(n) & = &\displaystyle{ \left(2^{j_n^{m_1(l')}-j_n^{m_1(l)}}\right)^{(1-\frac{2}{p})d} \int  \tilde \chi_l (y)
\ \tilde \chi_{l'}^{(\frac{p}{2}-1)} (x_n^{l,l'}(y))\ dy} \\\\
 & = &\displaystyle{ \left(2^{j_n^{m_1(l)}-j_n^{m_1(l')}}\right)^{(\frac{2}{p})d} \int  \tilde \chi_l (x_n^{l',l}(z))
\ \tilde \chi_{l'}^{(\frac{p}{2}-1)} (z)\ dz}\ .
\end{array}
\end{equation}
By orthogonality of $\tau_{m_1(l),n}$ and $\tau_{m_1(l'),n}$ (i.e. (\ref{six}), recall (\ref{m})), there are two possibilities:
\begin{enumerate}
\item One of the pre-factors in (\ref{ak}) tends to zero in which case $I^{l,l'}(n) \to 0$ as $n\to \infty$ since the functions in the integrands are bounded with compact support.
\item Both pre-factors are bounded, hence
$$2^{j_n^{m_1(l')}-j_n^{m_1(l)}} \in (c_1,c_2) \qquad \forall \ n$$
for some $c_2>c_1 >0$, and
$$\left|2^{j_n^{m_1(l')}-j_n^{m_1(l)}}k_n^{m_1(l)} - k_n^{m_1(l')}\right| \xto +\infty\ .$$
In this case, in light of (\ref{zzz}) we see that $I^{l,l'}(n) \equiv 0$ for sufficiently large $n$, since the supports of the two functions in the integrand in (\ref{ak}) become disjoint.
\end{enumerate}
Summing over a finite number of $l$ and $l'$ concludes the proof of (\ref{ag}) and of the theorem. \hfill $\Box$
\\\\
Note that we did not need to use the compact support of the $\psi^{(i)}$ (as was emphasized in \cite{jaffard}), only the compact support of the functions $\chi_\d$ (in fact, $\chi_\d \in \lprd$ would have been sufficient) and the relation (\ref{e}).  Hence, any well-localized basis functions $\psi^{(i)}$ would be sufficient for our purposes (e.g., the divergence-free wavelet bases of \cite{bf} or \cite{lr}).
\section{Other applications of the method}
\noindent
We can similarly (and in fact, more easily) deduce the following as well:
\\\\
Fix any $p \in [-\infty,+\infty]\backslash \{0\}$, and for any $a\geq 1$ set $s_{p,a}:= \frac{d}{a} - \frac{d}{p}$.  Bernstein's inequalities (\ref{aac}) immediately give us $\dot B^{s_{p,a}}_{a,q} \hookrightarrow \dot B^{s_{p,b}}_{b,r}$ whenever $1\leq a < b\leq +\infty$ and $1\leq q \leq r \leq +\infty$.  (The choice of $p$ dictates the common scaling of these spaces, which coincides with that of $L^p(\mathbb{R}^d)$ when $p>0$.)  If, moreover, we have the slight improvement that $r\geq \frac{b}{a} q$, the methods used to prove Theorem \ref{thm:a} easily enable us to prove the following similar result:
\begin{thm}\label{thm:aa}
Let $\{u_n\}_{n=1}^\infty$ be a bounded sequence in $\besa$.  There exists a sequence $\{\phi_l\}_{l=1}^\infty \subset \besa$ and for each $l\in \N$ a sequence $\{(j_n^l,k_n^l)\}_{n=1}^\infty \subset \Z \times \Z^d$, both depending on $\{u_n\}$, such that, after possibly passing to a subsequence in $n$,
\begin{equation}\label{profilesa}
u_n(x) = \sum_{l=1}^L (2^{j^l_n})^{{d/p}} \phi_l(2^{j^l_n}x - k^l_n) + r_n^L(x)
\end{equation}
for any $L \in \N$ where the following properties hold:
\begin{enumerate}
\item For $l\neq l'$, the sequences $\{(j_n^l,k_n^l)\}$ and $\{(j_n^{l'},k_n^{l'})\}$ are orthogonal in the following sense:
\begin{equation}\label{orth1a}
\left|\log \left(2^{(j_n^l - j_n^{l'})}\right)\right| + \left|2^{(j_n^l - j_n^{l'})}k_n^{l'} - k_n^l\right|  \xrightarrow[n\to\infty]{} +\infty
\end{equation}
\item The remainder $r_n^L$ satisfies the following smallness condition:
\begin{equation}\label{orth2a}
\lim_{L\to\infty} \left(\varlimsup_{n\to \infty} \|r_n^L\|_\besb \right) = 0
\end{equation}
\item There is a norm $\| \cdot \tilde \|_\besa$ which is equivalent to $\|\cdot \|_\besa$ such that for each $n \in \N$,
\begin{equation}\label{orth3a}
\left\| \left(\|(2^{j^l_n})^{{d/p}} \phi_l(2^{j^l_n}x - k^l_n) \tn_\besa\right)_{l=1}^\infty \right\|_{\ell^\tau} \leq \varliminf_{n'\to\infty} \|u_{n'}\tn_\besa
\end{equation}
where $\tau := \max \{a,q\}$, and for any $L \in \N$,
\begin{equation}\label{orth4a}
 \|r_n^L \tn_\besa \leq  \|u_{n}\tn_\besa + \circ(1) \quad \textrm{as} \quad n\to\infty\ .
\end{equation}
\end{enumerate}
\end{thm}
\noindent
The proof of Theorem \ref{thm:a} depended in part on the ``improved Sobolev-type'' embedding inequality\footnote{We believe moreover that in many settings of this type one can complete the proofs without explicitly using such an improved embedding inequality.  This will be addressed soon in \cite{ck}.} (\ref{j}).  The slight improvement on the distance between $r$ and $q$ is to allow us the following inequality, proved in Proposition \ref{prop:a3}, which is the analogue of (\ref{j}):
\begin{equation}\label{eq:sobc}
\|u\|_\besb \leq \|u\|^\alpha_\besa \|u\|^{1-\alpha}_\binf \ .
\end{equation}
The proofs of statements (\ref{profilesa}), (\ref{orth1a}), (\ref{orth2a}) and (\ref{orth4a}) are similar to (and, in fact, simpler than) the proofs of the analogous statements in Theorem \ref{thm:a} so we will neglect the details (see also \cite{ck}).  We will just mention that inequality (\ref{eq:sobc}) will replace inequality (\ref{j}), and one can replace the inequality in Lemma \ref{lemma1} by (\ref{decay2}).  This takes care of the key ingredients of the improved Sobolev-type inequality and the decay of non-linear wavelet projections.
\\\\
As a technical point, we also briefly mention the following:  Whereas in the proof of Theorem \ref{thm:a} one often needed to use Fatou's Lemma twice, once for an $L^q$-type norm and once for an $\ell^q$-type norm (for example in the proof of Lemma \ref{lemma:a4}), in the case of Besov spaces one proceeds in the same way but now both are $\ell^q$-type norms due to (\ref{c}), and so one no longer has to be cautious about pointwise limits in $\R^d$.  In this case, the ``change of variables" which allows the use of Fatou's lemma comes from the fact that within each profile, one can ``shift the labels" on the coefficients in a way depending on $n$ so that when one takes the Besov norm the only dependence on $n$ is in the coefficients (and not on the order in which they are summed).
\\\\
We turn therefore to the proof of (\ref{orth3a}), the proof of which is necessarily different since one cannot use the orthogonality of the scales and shifts which appear in the profile decomposition because our norms do not involve a space variable, only wavelet coefficients.  This will be compensated by the natural orthogonalities in the wavelet basis in spaces of the type given by (\ref{c}).
\\\\
The proof is due to the following basic inequality:
\\\\
Suppose $\LL$ is decomposed into a disjoint union of sets $E_i$:  $\LL = \cup_{i=1}^\infty E_i$, and suppose $f= \sum f_\l \psi_\l \in \besa$.  Then one has
$$
\left\| \left( \left\| \sum_{\l \in E_i} f_\l \psi_\l \right\|_X \right)_{i=1}^\infty \right\|_{\ell^\tau} \leq \left\| \sum_{\l \in \LL} f_\l \psi_\l \right\|_X
$$
where $\|\cdot \|_X := \| \cdot \tn_\besa$ and $\tau:= \max \{a,q\}$.  The inequality is a simple consequence of properties of $\ell^r$ spaces and definition (\ref{c}).
\\\\
Letting $E_i$ correspond to the wavelet components of $u_n$ which contribute to the rescaling of the $i$'th profile $\phi_i$ at a certain step of the iteration process, and bounding the term on the right by the full expansion of $u_n$, taking limits and using Fatou's lemma gives the desired result.  We refer the interested reader to the upcoming work \cite{ck} for more details.
\appendix
\section{Appendix}
\noindent
In this appendix we collect some known propositions regarding wavelets and embeddings of Besov and Lebesgue spaces, as well as a lemma (Lemma \ref{lemma:a4}) which is necessary for the proof of our main theorem.
\begin{prop}\label{prop:a1}
Let $\{a_\l\}$ correspond to some $f\in \lprd$, and for $A>0$ denote $\LL_A := \{ \l \in \LL \ | \ |a_\l| \geq A\ \}$.  Then $\#(\LL_A) < \infty$ for any $A>0$.  (In particular, $\exists \l_0 \in \LL$ such that $\sup_{\l \in \LL} |a_\l| = |a_{\l_0}|$.)
\end{prop}
\noindent
(For a countable set $\Omega$ we are denoting by $\#(\Omega)$ the number of elements in $\Omega$.)
\\\\
\textbf{Proof:}
\\\\
Fix $A>0$, and suppose to the contrary that $\#(\LL_A) = +\infty$.
\\\\
Suppose that there exists $\{\l_m = (i_m,j_m,k_m) \} \subset \LL_A$ such that $j_m \to -\infty$ as $m \to\infty$.  Then $|\d_m| = |\d(\l_m)| \to +\infty$ (recall that $\d(\l)$ is the cube associated to $j$ and $k$ if $\l = (i,j,k)$ and has volume $2^{-dj}$) and we can pass to a subsequence so that $|\d(\l_m)|$ increases so quickly that
\begin{equation}\label{aaa}
|\d^{(m)}| \geq \tfrac{1}{2} |\d(\l_m)| = \tfrac{1}{2}\cdot 2^{-dj_m}\ ,
\end{equation}
where we define $\d^{(m)} := \d(\l_m) \backslash \cup_{n=1}^{m-1} \d(\l_n)$.  Then for any $N\geq 2$,
\begin{equation}\label{aab}
\begin{array}{c}
\displaystyle{\|f\tn_\lprd^p   \geq \int_\rd \left\{ \sum_{m=1}^\infty |a_{\l_m}|^2 2^{\frac{2d}{p}j_m} \chi_{\d_m} \right\}^\frac{p}{2} \geq
\sum_{M=2}^{N+1} \int_{\d^{(m)}} \left\{ \sum_{m=1}^\infty |a_{\l_m}|^2 2^{\frac{2d}{p}j_m} \chi_{\d_m} \right\}^\frac{p}{2}} \\\\
\displaystyle{ \geq
\sum_{M=2}^{N+1} \int_{\d^{(m)}} \left\{ |a_{\l_M}|^2 2^{\frac{2d}{p}j_M} \chi_{\d_M} \right\}^\frac{p}{2} \geq \sum_{M=2}^{N+1} \frac{A^p}{2} = N \frac{A^p}{2}\ ,  }
\end{array}
\end{equation}
which implies (letting $N \to \infty$) that $\|f\|_\lprd = +\infty$ contrary to assumption.
\\\\
Similarly, if there exists $\{\l_m  \} \subset \LL_A$ such that $j_m \to +\infty$, for fixed $N \geq 2$ we can pass to a subsequence (depending on $N$) so that $|\d(\l_m)| \to 0$ rapidly enough that (\ref{aaa}) holds for $2\leq m \leq N+1$ after reversing the order of the first $N+2$ elements of the sequence.  By (\ref{aab}) therefore we have $\|f\tn_\lprd^p \geq N \frac{A^p}{2}$ which is false for a suitably large $N$ depending on $f$.
\\\\
Hence $i$ and $j$ can take only finitely-many values for $\l \in \LL_A$, so $\#(\LL_A) =+\infty$ implies that there exists $\{\l_m\} \subset \LL_A$ such that $|k_m| \to +\infty$.  Passing to a subsequence so that $j_m \equiv \bar{j} = \textrm{const}$ and $|k_m| \to +\infty$ sufficiently rapidly we can ensure that the $\d_m$ are mutually disjoint, so $\d^{(m)} = \d_m$ for all $m$.  Equations (\ref{aaa}) and (\ref{aab}) then again contradict the assumption that $f\in \lprd$ and the proposition is proved. \hfill $\Box$
\\\\
For the following two propositions, we will use the Littlewood-Paley description of Besov and Triebel-Lizorkin spaces in terms of the frequency localization operators $\D_j$, where $\widehat{\D_j f}$ is supported in a neighborhood of $\{|\xi| = 2^j\}$ for $j\in \Z$.  For example, one defines the Besov spaces by $\|f\|_{\dot B^s_{r,q}} = \left\| 2^{sj}\|\D_j f\|_{L^r(\rd)} \right\|_{\ell^q(\Z)}$.
\begin{prop}\label{prop:a2}
Fix $p,q,r \in \R$ such that $2\leq p \leq q,r \leq \infty$, and set ${s_{p,r}:= d(\tfrac{1}{r} - \tfrac{1}{p})}$ (note $s_{p,r} < 0$).  Then
$$\lprd \hookrightarrow \brq\ .$$
\end{prop}
\pagebreak
\noindent
\textbf{Proof:}
\\\\
Note first by Bernstein's inequalities that
\begin{equation}\label{aac}
2^{jd(\frac{1}{r}-\frac{1}{p})}\|\D_jf\|_{L^r(\rd)} \lesssim \|\D_jf\|_\lprd
\end{equation}
for
$$1\leq p \leq r \leq +\infty\ .$$
Taking $\ell^q$ norms of both sides of (\ref{aac}) gives $\dot B^0_{p,q} \hookrightarrow \bfin$.
\\\\
Note now that $\|f\|_\lprd \simeq \|f\|_{F^0_{p,2}} = \|\{\sum_j |\D_jf|^2\}^\frac{1}{2}\|_\lprd$.  Using the simple fact that $\|c_j\|_{\ell^\infty} \leq \|c_j\|_{\ell^\alpha}$ for any $\alpha >0$ implies\footnote{Indeed, $a>0$ and $\frac{a}{b} \in [0,1]$ gives
$\|c_j\|_{\ell^b} = \||c_j|^{1-\frac{a}{b}}|c_j|^\frac{a}{b}\|_{\ell^b} \leq  \||c_j|^{1-\frac{a}{b}}\|_{l^\infty}\||c_j|^\frac{a}{b}\|_{\ell^b}
\leq  \||c_j|^{1-\frac{a}{b}}\|_{l^{\frac{a}{1-\frac{a}{b}}}}\||c_j|^\frac{a}{b}\|_{\ell^b} = \|c_j\|_{\ell^a}$.} that $\|c_j\|_{\ell^b} \leq   \|c_j\|_{\ell^a}$ for $0< a \leq b \leq \infty$, we have
$$  \Big\{\sum_j |\D_jf|^p\Big\}^\frac{1}{p}(x) \leq \Big\{ \sum_j |\D_jf|^2\Big\}^\frac{1}{2}(x)$$
for any $x\in \rd$, as long as $p\geq 2$.  Taking $\lprd$ norms of both sides, Fubini's theorem now gives
\begin{equation}\label{a1}
\|f\|_{\dot B^0_{p,p}}= \left\|\Big\{\sum_j |\D_jf|^p\Big\}^\frac{1}{p}\right\|_\lp \leq \|f\|_{F_{p,2}^{0}} \lesssim \|f\|_\lp\ .
\end{equation}
Since $0<p\leq q \leq \infty$, we see therefore that
$$\lprd \hookrightarrow \dot B^0_{p,p} \hookrightarrow \dot B^0_{p,q} \hookrightarrow \bfin$$
and the proposition is proved. \hfill $\Box$
\begin{prop}\label{prop:a3}
Suppose $2\leq p < q,r \leq \infty$.  Then for any $\alpha \in (\max \{\frac{p}{r},\frac{p}{q}\},1)$, there exists $C=C(d,p,q,r,\alpha)$ such that
$$\|f\|_\brq \leq C \|f\|^\alpha_\lprd \|f\|^{1-\alpha}_\binf\ .$$
\end{prop}
\noindent
\textbf{Proof:}
\\\\
Fix such an $\alpha$ and set $a:=r\alpha$, $b:=q\alpha$.  Since
$$\frac{p}{r} < \alpha < 1 \qquad \Longrightarrow \qquad 1<a<r\ ,$$
H\" older's inequality gives
$$2^{jd(\frac{1}{r} - \frac{1}{p})} \|\D_jf\|_r \leq \left[2^{jd(\frac{1}{a} - \frac{1}{p})} \|\D_jf\|_a\right]^\alpha \left[ 2^{-j\frac{d}{p}} \|\D_jf\|_\infty\right]^{1-\alpha} \ .$$
Writing this as $\e_{r,j} \leq \e_{a,j}^\alpha \e_{\infty,j}^{1-\alpha}$, we have
$$\|\e_{r,j}\|_{\ell^q} \leq \|\e_{a,j}^\alpha \|_{\ell^q} \|\e_{\infty,j}^{1-\alpha}\|_{l^\infty} = \|\e_{a,j}\|_{\ell^b}^\alpha  \|\e_{\infty,j}\|_{l^\infty}^{1-\alpha}$$
which means that
$$\|f\|_\bfin \leq \|f\|_{\dot B^{s_{p,a}}_{a,b}}^\alpha \|f\|^{1-\alpha}_\binf \ .$$
But by Proposition \ref{prop:a2}, $\|f\|_{\dot B^{s_{p,a}}_{a,b}} \lesssim \|f\|_\lprd$ provided only that $2\leq p \leq a,b \leq +\infty$.  This is exactly guaranteed by our assumptions on $\alpha$, and the proposition follows. \hfill $\Box$
\\\\
Finally, we prove the following crucial lemma:
\begin{lemma}\label{lemma:a4}
For each $n$, define $\phi_n \in L^p(\R^d)$ by
$$\phi_n = \sum_{\mu \in \N} a_n^\mu \psi_{\l_n^\mu}$$ and suppose the following are satisfied:
\begin{enumerate}
\item There exists $C>0$ such that $\|\phi_n\|_p \leq C$ for all $n \in \N$,
\item for each $\mu\in \N$, there exists $a^\mu \in \R$ such that $a_n^\mu \to a^\mu$ as $n\to \infty$, and
\item for each $\mu\in \N$, there exists $j^{(\mu,1)} \in \Z$ and $b^{(\mu,1)} \in \R^d$ such that
$$j_n^\mu - j_n^1 \equiv j^{(\mu,1)} \qquad \textrm{and} \qquad k_n^\mu - 2^{j^{(\mu,1)}} k_n^1 \equiv b^{(\mu,1)} \qquad \textrm{for all}\ n\in \N\ ,$$
and $\l_n^\mu = (i_\mu,j_n^\mu, k_n^\mu)$ for some $i_\mu \in \{1,\dots,2^d-1\}$ for all $\mu,n\in\N$.
\end{enumerate}
Then $\tilde\phi_n := \sum_{\mu\in \N} a^\mu \psi_{\l_n^\mu} \in L^p(\R^d)$ for any $n\in \N$ and there exists $\phi \in L^p(\R^d)$ such that $\tilde\phi_n \equiv \tau_{n,1}\phi$ for all $n\in \N$.
\end{lemma}
\ \\
\noindent
\textbf{Proof:}
\\\\
Note that
$$
\begin{array}{rcl}
\d_n^\mu & = & \{2^{j_n^\mu}x - k_n^\mu \in [0,1)^d\ \} \\\\
& = & \{2^{j_n^1 + j^{(\mu,1)}}x - 2^{j^{(\mu,1)}}k_n^1 - b^{(\mu,1)} \in [0,1)^d\ \} \\\\
& = & \{2^{j^{(\mu,1)}}( 2^{j_n^1}x -k_n^1) - b^{(\mu,1)} \in [0,1)^d\ \} \\\\
& = & \{\tau^{(\mu,1)}(\tau_{n,1}(x)) \in [0,1)^d\ \}\\\\
& = & (\tau_{n,1})^{-1}[ (\tau^{(\mu,1)})^{-1}([0,1)^d)]
\end{array}
$$
where $\tau^{(\mu,1)}(x) := 2^{j^{(\mu,1)}}x - b^{(\mu,1)}$.  Defining $\d^{(\mu,1)} := (\tau^{(\mu,1)})^{-1}([0,1)^d)$, we can now easily see that
$$\chi_{\d_n^\mu}(x) = \chi_{\d^{(\mu,1)}}(\tau_{n,1}(x))\ .$$
\pagebreak
Using now the change of variables $y = \tau_{n,1}(x) = 2^{j_n^1}x - k_n^1$, we see that
$$I_n := \int \left\{ \sum_\mu |a_n^\mu|^2 2^{\frac{2d}{p}j_n^\mu} \chi_{\d_n^\mu}(x)   \right\}^{{p/2}} dx =
\int \left\{ \sum_\mu |a_n^\mu|^2 2^{\frac{2d}{p}j_n^\mu} \chi_{\d^{(\mu,1)}}(\tau_{n,1}(x))   \right\}^{{p/2}} dx =
$$
$$=\int \left\{ \sum_\mu |a_n^\mu|^2 2^{\frac{2d}{p}j^{(\mu,1)}} \chi_{\d^{(\mu,1)}}(y)   \right\}^{{p/2}} dy
$$
Note now that the only dependence on $n$ in the last integral is in the terms $a_n^\mu$.
\\\\
Fix first any $y\in \R^d$ such that the following is finite for every $n$ (it is true almost everywhere because $\phi_n \in L^p$ implies $I_n$ is finite for each $n$, and a countable union of sets of measure zero has measure zero), and write
$$\sum_\mu |a_n^\mu|^2 2^{\frac{2d}{p}j^{(\mu,1)}} \chi_{\d^{(\mu,1)}}(y) = \int_{\N} |a_n^\mu|^2 2^{\frac{2d}{p}j^{(\mu,1)}} \chi_{\d^{(\mu,1)}}(y) d\alpha (\mu)$$
where $\alpha$ is the discrete measure.  By Fatou's Lemma, we have
$$\int_{\N} \left\{\liminf_{n\to\infty} |a_n^\mu|^2 2^{\frac{2d}{p}j^{(\mu,1)}} \chi_{\d^{(\mu,1)}}(y)\right\} d\alpha (\mu) \leq \liminf_{n\to\infty} \int_{\N} |a_n^\mu|^2 2^{\frac{2d}{p}j^{(\mu,1)}} \chi_{\d^{(\mu,1)}}(y) d\alpha (\mu)\ ,$$
hence
$$\left\{ \sum_m |a^\mu|^2 2^{\frac{2d}{p}j^{(\mu,1)}} \chi_{\d^{(\mu,1)}}(y)\right\}^\frac{p}{2} \leq \left\{ \liminf_{n\to\infty} \sum |a_n^\mu|^2 2^{\frac{2d}{p}j^{(\mu,1)}} \chi_{\d^{(\mu,1)}}(y)\right\}^\frac{p}{2} =$$$$=  \liminf_{n\to\infty} \left\{ \sum |a_n^\mu|^2 2^{\frac{2d}{p}j^{(\mu,1)}} \chi_{\d^{(\mu,1)}}(y)\right\}^\frac{p}{2} \ .$$
Integrating now on both sides of the previous inequality in $y$ and applying Fatou's lemma again with respect to the usual Lebesgue measure on $\R^d$, we have
$$I := \int \left\{ \sum_\mu |a^\mu|^2 2^{\frac{2d}{p}j^{(\mu,1)}} \chi_{\d^{(\mu,1)}}(y)   \right\}^{{p/2}} dy \leq \liminf_{n\to\infty} I_n$$
and we know by equivalence of the wavelet norms and the usual Lebesgue norms that $I_n \lesssim \|\phi_n\|_p^p \leq C$ hence ${I \lesssim C <\infty}$.  Changing back to the variable $x$ in $I$ (which gives the same formula as the original definition of $I_n$ as an integral in $x$ only with $a_n^\mu$ replaced by $a^\mu$), we see again by the equivalence of norms that
$\|\tilde{\phi}_n\|_p^p \lesssim I < \infty$.  (If one is uncomfortable considering the wavelet expansion of a potentially non-existence element of $\lprd$, one can alternatively show that the relevant sequences of partial sums are Cauchy by similar calculations to obtain the existence of $\tilde \phi_n$.)
\\\\
It is obvious from the setup that $\tau_{n,1}^{-1} \tilde{\phi}_n$ is independent of $n$, hence the last statement in the claim.  Indeed, we must have
$$\tau_{n,1}^{-1} \tilde{\phi}_n = \sum_\mu a^\mu \tau^{(\mu,1)}\psi^{(i_\mu)}$$
and the lemma is proved.
\hfill $\Box$

\bibliography{bibliography}
\bibliographystyle{plain}
\end{document}